\documentclass[a4j]{article}
\usepackage{amssymb}
\usepackage{amsmath}

\newtheorem{teiri}{Theorem}
\newtheorem{prop}{Proposition}

\title{The Hannan-Quinn Proposition for Linear Regression.}
\author{Joe Suzuki\thanks{joe@suzuki.email.ne.jp}\\
Department of Mathematics, Osaka University\thanks{Toyonaka, Osaka 560-0043, Japan.}}
\date{MSC2010: 62J05, 62M10}
\begin{document}

\large 

\maketitle

\section*{Abstract}
We consider the variable selection problem in linear regression.
Suppose that we have a set of random variables $X_1,\cdots,X_m,Y,\epsilon$
such that $Y=\sum_{k\in \pi}\alpha_kX_k+\epsilon$ with $\pi\subseteq \{1,\cdots,m\}$ and
$\alpha_k\in {\mathbb R}$ unknown, and $\epsilon$ is independent of any linear combination of $X_1,\cdots,X_m$.
Given actually emitted $n$ examples $\{(x_{i,1}\cdots,x_{i,m},y_i)\}_{i=1}^n$ emitted from $(X_1,\cdots,X_m, Y)$,
we wish to estimate the true $\pi$ using information criteria in the form of $H+(k/2)d_n$,
where $H$ is the likelihood with respect to $\pi$ multiplied by $-1$, and $\{d_n\}$ is a positive real sequence.
If $d_n$ is too small, we cannot obtain consistency because of overestimation.
For autoregression, Hannan-Quinn proved that, in their setting of $H$ and $k$, 
the rate $d_n=2\log\log n$ is the minimum  satisfying strong consistency.
This paper solves the statement affirmative for linear regression as well which has a completely different setting.

\section*{Keywords}
Hannan-Quinn, linear regression, the law of iterated logarithms, strong consistency, information criteria, model selection.

\section{Introduction}

Information criteria such as AIC, MDL/BIC are used for problems in model selection, and
each problem is associated with estimating how many independent parameters exist from given finite examples:
on how many variables another variable depends in linear regression (LR);
on how many previous variables the subsequent variable depends on in auto regression (AR),
etc. 

For each model $g$, we evaluate two factors:
\begin{enumerate}
\item How well the examples explain the model $g$; and
\item How simple the model $g$ is.
\end{enumerate}
and balance them numerically.
Let $\{d_n\}_{n=1}^\infty$ be nonnegative reals such that
$d_n/n \rightarrow 0$,
$H(g)$ 
the empirical entropy which is the maximum likelihood multiplied by $(-1)$,
and $k(g)$ the number of parameters in model $g$.
By information criteria, we mean the quantity
\begin{equation}\label{eq898}
H(g)+\frac{k(g)}{2}d_n\ ,
\end{equation}
and we estimate the model $g$ by finding one with the minimum value.
For example, $d_n=2$ for AIC, and $d_n=\log n$ for MDL/BIC.
Hence, information criteria exist as many as sequences $\{d_n\}_{n=1}^\infty$, so 
it is impossible to list all of information criteria in the form of (\ref{eq898}).

In model selection, in particular for theoretical analyses, we often discuss if
consistency holds for each $\{d_n\}$,
namely, if a sequence of selected models converges to the correct one as
$n\rightarrow \infty$ in the following senses:
\begin{enumerate}
\item the probability of the selected model for each $n$ being correct converges to one (weakly consistent), and
\item the set (event) of infinite sequences in which at most a finite number of errors
occur has probability one (strongly consistent).
\end{enumerate}
Although both properties are satisfied in MDL/BIC ($d_n=\log n$),
however, none of the two  are satisfied in AIC ($d_n=2$).
In general, if 
$d_n$ is too small, strong consistency is not obtained because of overestimation.

This paper addresses the minimum order of  $\{d_n\}$ satisfying strong consistency although
seeking such a condition is of theoretical interest in model selection (in fact, many information criteria
are to be satisfactory even if consistency is not achieved). 

The definitions of empirical entropy and the number of parameters are different
in each problem to be considered.
In 1979, Hannan-Quinn proved that for AR
$d_n=2\log\log n$ is the minimum order satisfying strong consistency (Hannan-Quinn proposition).
However, the same $d_n=2\log\log n$ has been applied to other problems as well as AR.
In fact, 
the proof of the Hannan-Quinn proposition essentially depends on 
the properties of the AR problem, which is clear from the original paper by Hannan-Quinn, and
the Hannan-Quinn proposition was not proved for any other problem including the LR problem.
On the contrary, without noticing such a matter, the information criterion HQ was applied to those problems.

Recently, the Hannan-Quinn proposition has been proved for estimating classification rules which has many applications such as
Markov order estimation, data mining, pattern recognition (Suzuki, 2006).

This paper shows that the Hannan-Quinn proposition is true for estimating
dependencies in LR, which seems to be of great significance.
Otherwise, there would be no reason to use HQ in LR.
Several authors suggested that $d_n=c\log\log n$ with some positive
 constant $c$ would be enough (Rao-Wu, 1989). So, there has been evidence
 that the proposition is true although no formal proof appeared.
This paper proves that such a $c$ is any constant strictly greater than two.

In Section 2, we briefly overview how the Hannan-Quinn proposition was proved in 
AR.
In Section 3, we derive the asymptotic error probability of model selection in LR
when information criteria are applied, which will be an important step to prove the main result.
In Section 4, we give a proof of 
the Hannan-Quinn proposition for LR.
Section 5 summarizes the results in this paper and gives a future problem.

Throughout the paper, we denote by $X(\Omega)$ the image $\{X(\omega)|\omega\in \Omega\}$ of a random variable $X: \Omega\rightarrow {\mathbb R}$,
where $\Omega$ is the underlying sample space.

\section{Auto Regression}


Let 
$\{W_i\}_{i=-\infty}^\infty$ be a sequence of independent and identically distributed 
random variables with expectation zero and variance one, and let 
$\{X_i\}_{i=-\infty}^\infty$ be defined by
$$X_i=\sum_{j=1}^k\lambda_jX_{i-j}+W_i$$
and a nonnegative real sequence $\{\lambda_i\}_{i=1}^k$,
where we assume the expectation of each $X_i$ to be  zero.
Since $\{X_i\}$ is stationary, we obtain for $m\geq 0$, the following equation (Yule-Waker)
$$\gamma_{m}=\sum_{j=1}^k\lambda_j \gamma_{m-j}+\delta_{0m}\sigma^2_k\ ,$$
where $\gamma_m:=E X_iX_{i+m}$ does not depend on $i$.
Using Cramer's formula, and from the values of $\{\gamma_m\}_{m=0}^k$, we 
obtain the values of $\lambda_0:=\sigma_k^2$ and $\{\lambda_m\}_{m=1}^k$ as a solution of
the $(k+1)\times (k+1)$ linear equations:
$$
\left[
\begin{array}{ccccc}
-1&{\gamma}_1 &{\gamma}_2&\cdots &{\gamma}_k\\
0 &{\gamma}_0 &{\gamma}_1&\cdots &{\gamma}_{k-1}\\
0 &{\gamma}_1 &{\gamma}_0&\cdots &{\gamma}_{k-2}\\
\vdots &\vdots &\vdots&\vdots &\vdots\\
0 &{\gamma}_{k-1} &{\gamma}_{k-2}&\cdots &{\gamma}_{0}\\
\end{array}
\right]
\left[
\begin{array}{c}
{\sigma}_k^2\\
{\lambda}_{1,k}\\
{\lambda}_{2,k}\\
\vdots\\
{\lambda}_{k,k}\\
\end{array}
\right]=
\left[
\begin{array}{c}
-{\gamma}_0\\
-{\gamma}_{1}\\
-{\gamma}_{2}\\
\vdots\\
-{\gamma}_{k}\\
\end{array}
\right]\ .
$$

Since the values of $\{\gamma_m\}_{m=0}^k$ are generally unknown,
we need to estimate 
$$\bar{x}:=\frac{1}{n}\sum_{i=1}^nx_i$$
and
$$\hat{\gamma}_m:=\hat{\gamma}_{-m}:=\frac{1}{n}\sum_{i=1}^{n-m}(x_i-\bar{x})(x_{i+m}-\bar{x})$$
from the examples
$$x^n=(x_1,\cdots,x_n)\in X_1(\Omega)\times \cdots \times X_n(\Omega)\ .$$
Then, we obtain the Yule-Walker equation as follows:
\begin{equation}\label{eq1-01}
\left[
\begin{array}{ccccc}
-1&\hat{\gamma}_1 &\hat{\gamma}_2&\cdots &\hat{\gamma}_k\\
0 &\hat{\gamma}_0 &\hat{\gamma}_1&\cdots &\hat{\gamma}_{k-1}\\
0 &\hat{\gamma}_1 &\hat{\gamma}_0&\cdots &\hat{\gamma}_{k-2}\\
\vdots &\vdots &\vdots&\vdots &\vdots\\
0 &\hat{\gamma}_{k-1} &\hat{\gamma}_{k-2}&\cdots &\hat{\gamma}_{0}\\
\end{array}
\right]
\left[
\begin{array}{c}
\hat{\sigma}_k^2\\
\hat{\lambda}_{1,k}\\
\hat{\lambda}_{2,k}\\
\vdots\\
\hat{\lambda}_{k,k}\\
\end{array}
\right]=
\left[
\begin{array}{c}
-\hat{\gamma}_0\\
-\hat{\gamma}_{1}\\
-\hat{\gamma}_{2}\\
\vdots\\
-\hat{\gamma}_{k}\\
\end{array}
\right]\ .
\end{equation}
In particular, if the order $k$ is unknown, we solve the above linear equation for each $k$ to calculate
the value of
\begin{equation}\label{eq425}
L(x^n,k)=\frac{n}{2}\log \hat{\sigma}_k^2 + \frac{k}{2}d_n\ .
\end{equation}
We estimate the true $k=k^*$ by the one $k=\hat{k}$ that minimizes (\ref{eq425}). 
This process is called estimating the AR order.
Then, we also obtain the solutions  $\hat{\lambda}_{0,\hat{k}}:=\hat{\sigma}^2_{\hat{k}}$
and $\{\lambda_{m,{\hat{k}}}\}_{m=1}^{\hat{k}}$ of (\ref{eq1-01}) with $k=\hat{k}$.

In general, 
$$\hat{\sigma}_k^2=\{1-\hat{\lambda}_{k,k}^2\}\hat{\sigma}_{k-1}^2\ ,$$
thus for 
each $k=1,2,\cdots$, we have
\begin{eqnarray}
&&2\{L(x^n,k)-L(x^n,k-1)\}\nonumber\\&=&{n}\log \frac{\hat{\sigma}_k^2}{\hat{\sigma}_{k-1}^2} + d_n\nonumber\\
&\leq& -n(1-\frac{\hat{\sigma}_{k}^2}{\hat{\sigma}_{k-1}^2})+d_n\label{eq31}\\
&=&-n\hat{\lambda}_{k,k}^2+d_n\ . \nonumber
\end{eqnarray}
As $n\rightarrow \infty$, for 
$k\leq  k^*$, 
$\displaystyle \frac{\hat{\sigma}_{k}^2}{\hat{\sigma}_{k-1}^2}$
almost surely converges to a value less than one.
Thus, from (\ref{eq31}), we have with probability one
$$L(x^n,0)>L(x^n,1)>\cdots >L(x^n,k^*-1)>L(x^n,k^*)\ .$$
On the other hand, for 
$k\geq k^*+1$, 
$\displaystyle \frac{\hat{\sigma}_{k}^2}{\hat{\sigma}_{k-1}^2}$
almost surely converges to one.
Hannan-Quinn(1979) proved from the law of iterated logarithms that
$$\frac{\hat{\lambda}_{k,k}^2}{2n^{-1}\log\log n}\leq 1$$
with probability one, and that for $d_n=2c\log\log n$ ($c>1$),
$$L(x^n,k^*)<L(x^n,k^*+1)<\cdots $$
with probability one.

\section{Linear Regression}

Let $X_1,\cdots,X_m$ be random variables such that there are no linear relations:
any linear combination of $X_1,\cdots,X_m$ cannot be zero with probability one.
Let $\epsilon\sim {\cal N}(0,\sigma^2)$  be a normal random variable with expectation zero and variance $\sigma^2>0$,
and 
$$Y:=\sum_{j=1}^p\alpha_jX_j+\epsilon\ ,$$
where ${\boldsymbol \alpha}:=[\alpha_1,\cdots,\alpha_p]^T\in {\mathbb R}^p$ ($0\leq p\leq m$).
We assume that 
$\epsilon$ is independent of any linear combination of $X_1,\cdots,X_m$.

Suppose we do not know the values of order $p$ and coefficients $\boldsymbol \alpha$,
and that we are given independently emitted 
$n$ examples 
$$z^n:=\{[y_i,x_{i,1},\cdots,x_{i,m}]\}_{i=1}^n$$ with
$$y_i\in Y(\Omega), [x_{i,1},\cdots,x_{i,m}]\in X_1(\Omega)\times \cdots \times X_m(\Omega)\ ,$$
where $\{[x_{1,j},\cdots,x_{n,j}]\}_{j=1}^m$ are to be linearly independent.
If we define 
$${\boldsymbol X}_p:=\left[
\begin{array}{ccc}
x_{1,1}&\ldots&x_{1,p}\\
\vdots&\ddots&\vdots\\
x_{n,1}&\ldots&x_{n,p}
\end{array}
\right]
,\ 
{\bold y}:=
\left[
\begin{array}{c}
y_1\\
\vdots\\
y_n
\end{array}
\right]
,\ 
{\boldsymbol \epsilon}:=\left[
\begin{array}{c}
\epsilon_1\\
\vdots\\
\epsilon_n
\end{array}
\right]\ ,
$$
we can write 
${\boldsymbol y}={\boldsymbol X}_p{\boldsymbol \alpha}+{\boldsymbol \epsilon}$.
Suppose that we estimate $p$ by $q$ ($0\leq q\leq m$). If 
we wish to minimize the quantity $\sum_{i=1}^n (y_i-\sum_{j=1}^q\hat{\alpha}_{jq}x_{ij})^2$ given the $n$ examples, then
$\hat{\boldsymbol \alpha}_q=[\hat{\alpha}_{1,q},\cdots,\hat{\alpha}_{q,q}]^T:=({\boldsymbol X}_q^T{\boldsymbol X}_q)^{-1}{\boldsymbol X}_q^T{\boldsymbol y}$
is the exact solution (minimum square error estimation), where
$${\boldsymbol X}_q:=\left[
\begin{array}{ccc}
x_{1,1}&\ldots&x_{1,q}\\
\vdots&\ddots&\vdots\\
x_{n,1}&\ldots&x_{n,q}
\end{array}
\right]
$$

\subsection{Idempotent Matrices}
Suppose $p\leq q$.
 If we define  $P_q:={\boldsymbol X}_q({\boldsymbol X}_q^T{\boldsymbol X}_q)^{-1}{\boldsymbol X}_q^T$, we have
$$P_q^2=P_q$$
and
$$(I-P_q)^2=I-P_q\ ,$$
so that the square error is expressed by
\begin{eqnarray*}
S_q&:=&\sum_{i=1}^n (y_i-\sum_{j=1}^q\hat{\alpha}_{j,q}x_{i,j})^2\\
&=&||{\boldsymbol y}-{\boldsymbol X}_q\hat{\boldsymbol \alpha}_q||^2\\
&=&||(I-P_q){\boldsymbol y}||^2\\
&=&{\boldsymbol y}^T(I-P_q){\boldsymbol y}\ .
\end{eqnarray*}
Similarly, if $q=p$, for 
$P_p:={\boldsymbol X}_p({\boldsymbol X}_p^T{\boldsymbol X}_p)^{-1}{\boldsymbol X}_p^T$ and 
$\hat{\boldsymbol \alpha}_p =[\hat{\alpha}_{1,p},\cdots,\hat{\alpha}_{p,p}]^T:=
({\boldsymbol X}_p^T{\boldsymbol X}_p)^{-1}{\boldsymbol X}_p^T{\boldsymbol y}$,
the square error is expressed by
$$S_p={\boldsymbol y}^T(I-P_p){\boldsymbol y}\ .$$

Thus, the difference between the square errors is
\begin{equation}\label{eq789}
S_p-S_q={\boldsymbol y}^T(I-P_q){\boldsymbol y}-{\boldsymbol y}^T(I-P_q){\boldsymbol y}={\boldsymbol y}^T(P_q-P_p){\boldsymbol y}\ .
\end{equation}
On the other hand, we have
\begin{eqnarray*}
P_q^T&=&(X_q^T)^T\{(X_q^TX_q)^{-1}\}^TX_q^T\\
&=&X_q\{(X_q^TX_q)^T\}^{-1}X_q^T=P_q
\end{eqnarray*}
and $P_p^T=P_p$.
From $P_qX_p=X_p$, $P_pX_p=X_p$, we obtain
$$P_qP_p=P_qX_p(X_p^TX_p)^{-1}X_p^T=X_p(X_p^TX_p)^{-1}X_p^T=P_p$$
and
$$P_pP_q=P_p^TP_q^T=(P_qP_p)^T=P_p^T=P_p\ .$$
Thus, not just for 
$P_p, I-P_p$ but also for $P_q-P_p$, the property
$$(P_q-P_p)^2=P_q^2-P_qP_p-P_pP_q+P_p^2=P_q-P_p$$
holds. Such square matrices satisfying the property are called idempotent matrices (Chatterjee-Hadi, 1987).

In general, for idempotent matrix $P\in {\mathbb R}^{n\times n}$,
the inner product $(Px, (I-P)x)=0$ for any $x=Px+(I-P)x \in {\mathbb R}^n$, so that
the eigenspaces are
\begin{enumerate}
\item $V_1:=\{Px|x\in {\mathbb R}^n\}$ with ${\rm dim}(V_1)={\rm rank}(P)$, and 
\item $V_0:=\{(I-P)x|x\in {\mathbb R}^n\}$ with ${\rm dim}(V_0)=n-{\rm rank}(P)$.
\end{enumerate}
Since the eigenvalues are one and zero,
the multiplicity of eigenvalue one is the same as the trace.
Notice that for $(X_q^TX_q)=[y_{jk}]$ and $(X_q^TX_q)^{-1}=[z_{jk}]$,
$$trace(P_q)=trace(X_q(X_q^TX_q)^{-1}X_q^T)
=\sum_{i=1}^n\sum_{j=1}^q\sum_{k=1}^qx_{ij}z_{jk}x_{ki}
=\sum_{j=1}^q\sum_{k=1}^q y_{kj}z_{jk} =\sum_{k=1}^q 1=q\ ,$$
and $trace(P_p)=p$, so that
we have the following table.
\begin{center}
\begin{tabular}{c|c|c|c|c}
\hline
$P$&${\rm trace}(P)$&${\rm dim}(V_1)$&${\rm dim}(V_0)$&${\rm rank}(P)$\\
\hline
$P_p$&$p$&$p$&$n-p$&$p$\\
$I-P_p$&$n-p$&$n-p$&$p$&$n-p$\\
$P_q-P_p$&$q-p$&$q-p$&$n-q+p$&$q-p$\\
\hline
\end{tabular}
\end{center}


\subsection{Error probability in model selection}

\begin{prop} If $p<q$, 
$\displaystyle \frac{S_p-S_q}{S_p/n}$ asymptotically obeys the $\chi^2$ distribution with freedom $q-p$.
\end{prop}
Proof: 
Given ${\boldsymbol X}_p$, we choose 
an orthogonal matrix  $U=[{\boldsymbol u}_1,\cdots,{\boldsymbol u}_{n}]$ of $I-P_p$ so that
$U_1=<{\boldsymbol u}_1,\cdots,{\boldsymbol u}_{n-p}>$
and
$U_0=<{\boldsymbol u}_{n-p+1},\cdots,{\boldsymbol u}_{n}>$ are the eigenspaces of
eigenvalues one and zero, respectively.
Notice that
\begin{equation}\label{eq656}
(I-P_p){\boldsymbol y}=
{\boldsymbol y}-(X_p{\boldsymbol \alpha}+P_p{\boldsymbol \epsilon})
={\boldsymbol \epsilon}-P_p{\boldsymbol \epsilon}=(I-P_p){\boldsymbol \epsilon}\ .
\end{equation}
For $j=1,\cdots,n-p$, multiplying ${\boldsymbol u}_j^T$ in both hands from left, we get a normal random variable
$$z_j:={\boldsymbol u}_j^T{\boldsymbol y}={\boldsymbol u}_j^T{\boldsymbol \epsilon}\ .$$

Since the expectation and variance of $\epsilon_i$ are zero and $\sigma^2$ (independent), and
$$
{\boldsymbol u}^T_j{\boldsymbol u}_k
=\left\{
\begin{array}{ll}
1,&j=k\\
0,&j\not=k
\end{array}
\right. \ ,
$$
we have $E[z_j]=0$ and
$$E[z_jz_k]=E[{\boldsymbol u}_j^T{\boldsymbol \epsilon}\cdot {\boldsymbol u}_k^T{\boldsymbol \epsilon}]=\sigma^2
{\boldsymbol u}_j^T{\boldsymbol u}_k
=\left\{
\begin{array}{ll}
\sigma^2,&j=k\\
0,&j\not=k
\end{array}
\right. \ .
$$
Thus, from the strong law of large numbers, with probability one as $n\rightarrow \infty$,
\begin{equation}\label{eq66}
\frac{1}{n}S_p=\frac{1}{n}\sum_{j=1}^{n-p}z_j^2\rightarrow \sigma^2\ .
\end{equation}

On the other hand, given ${\boldsymbol X}_q$,
we choose an orthogonal matrix $V=[{\boldsymbol v}_1,\cdots,{\boldsymbol v}_{n}]$ of 
$P_q-P_p$ so that 
$V_1=<{\boldsymbol v}_1,\cdots,{\boldsymbol v}_{q-p}>$
and
$V_0=<{\boldsymbol v}_{q-p+1},\cdots,{\boldsymbol v}_{n}>$ 
 are the eigenspaces of eigenvalues one and zero, respectively.
Notice that from (\ref{eq656}), we have 
$$(P_q-P_p){\boldsymbol y}=P_q(I-P_p){\boldsymbol y}=P_q(I-P_p){\boldsymbol \epsilon}=(P_q-P_p){\boldsymbol \epsilon}\ .$$
For $j=1,\cdots,q-p$, multiplying ${\boldsymbol v}_j$ in both hands from left, we get a normal random variable
$$r_j:={\boldsymbol v}_j^T{\boldsymbol y}=
{\boldsymbol v}_j^T{\boldsymbol \epsilon}\ .$$

Since the expectation and variance of $\epsilon_i$ are zero and $\sigma^2$ (independent), and
$$
{\boldsymbol v}^T_j{\boldsymbol v}_k
=\left\{
\begin{array}{ll}
1,&j=k\\
0,&j\not=k
\end{array}
\right. \ ,
$$ 
we have $E[r_j]=0$ and
$$E[r_jr_k]=E[{\boldsymbol v}_j^T{\boldsymbol \epsilon}\cdot {\boldsymbol v}_k^T{\boldsymbol \epsilon}]=\sigma^2
{\boldsymbol v}_j^T{\boldsymbol v}_k^T
=\left\{
\begin{array}{ll}
\sigma^2,&j=k\\
0,&j\not=k
\end{array}
\right. \ .
$$

Hence, as $n\rightarrow \infty$,
\begin{equation}\label{eq67}
\frac{S_p-S_q}{\sigma^2}=\sum_{j=1}^{q-p}\frac{r_j^2}{\sigma^2}\sim \chi^2_{q}
\end{equation}
where the fact that the square sum of $q-p$ independent random variables with the standard normal distribution
 obeys the $\chi^2$ distribution of freedom $q-p$ has been applied.
Equations (\ref{eq66})(\ref{eq67}) imply Proposition 1.
\begin{flushright}(Q. E. D.)\end{flushright}

In the sequel, for $\pi \subseteq \{1,\cdots,m\}$, we write
the square error of 
$\{X_j\}_{j\in \pi}$ and $Y$ by $S(\pi)$, and put
$$L(z^n,\pi):=n\log S(\pi)+\frac{k(\pi)}{2}d_n$$ 
and $k(\pi)=|\pi|$, given $z^n=\{[y_i,x_{i,1},\cdots,x_{i,m}]\}_{i=1}^n$.
 Let 
$\pi_* \subseteq \{1,\cdots,m\}$ be the true $\pi$.
\begin{teiri}\rm
For $\pi \supset \pi_*$, 
the probability of $L(z^n,\pi)<L(z^n,\pi_*)$ is 
$$\int_{n\{1-\exp[- \frac{k(\pi)-k(\pi_*)}{2n}d_n]\}}^\infty f_{k(\pi)-k(\pi_*)}(x)dx\ ,$$
where $f_l$ is the probability density function of the $\chi^2$ distribution of freedom $l$.
\end{teiri}
Proof: Notice that
\begin{eqnarray}
&&2\{L(z^n,\pi)-L(z^n,\pi_*)\}\nonumber\\
&=&{2n}\log \frac{S(\pi)}{S(\pi_*)} + \{k(\pi)-k(\pi_*)\}d_n\nonumber
\\
&=& 2n\log(1-\frac{S(\pi_*)-S(\pi)}{S(\pi_*)})+\{k(\pi)-k(\pi_*)\}d_n \label{eq79}\ ,
\end{eqnarray}
so that
\begin{eqnarray}\label{eq886}
L(z^n,\pi)<L(z^n,\pi_*)&\Longleftrightarrow& \frac{S(\pi_*)-S(\pi)}{S(\pi_*)/n}>
n\{1-\exp[- \frac{k(\pi)-k(\pi_*)}{2n}d_n]\}
 \label{eq29}\ .
\end{eqnarray}
From Proposition 2, we obtain Theorem 1.
\begin{flushright}(Q. E. D.)\end{flushright}

Hereafter, we do not assume that $\epsilon_i \sim {\cal N}(0,\sigma^2)$ but that $\epsilon_i$
is an independently identically distributed random variable with expectation zero and variance $\sigma^2$.

\begin{teiri}\rm
For $\pi\not\supseteq \pi_*$, $L(x^n,\pi)>L(x^n,\pi_*)$
with probability one as 
$n\rightarrow \infty$.
\end{teiri}
Proof: 
Suppose $q<p$. Given ${\boldsymbol X}_p$, we choose an orthogonal matrix $W:=[{\boldsymbol w}_1,\cdots,{\boldsymbol w}_n]$
of $P_p-P_q$ so that $W_1=<{\boldsymbol w}_1,\cdots,{\boldsymbol w}_{p-q}>$ and $W_0=<{\boldsymbol w}_{p-q+1},\cdots,{\boldsymbol w}_{n}>$
are the eigenspaces of eigenvalue one and zero, respectively.
Since $\{\hat{\alpha}_{j,p}\}_{j=1}^p$ are 
strongly consistent estimators (Lai-Robbins-Wei, 1978), we have for $j=1,\cdots p-q$ with probability one as $n\rightarrow \infty$
\begin{eqnarray*}
s_j&:=&\sum_{i=1}^n w_{i,j}y_i=\sum_{i=1}^n w_{ij}\{\sum_{k=1}^p x_{ik}\hat{\alpha}_{k,p}+y_i-\sum_{k=1}^p x_{ik}\hat{\alpha}_{k,p}\}\\
&\rightarrow& \sum_{i=1}^n w_{ij}(\sum_{k=1}^p x_{ik}{\alpha}_{k}+\epsilon_i)\\
&\rightarrow& \sum_{i=1}^n w_{ij}(\sum_{k=q+1}^p x_{ik}{\alpha}_{k}+\epsilon_i)\ ,
\end{eqnarray*}
where ${\boldsymbol w}_j:=[w_{1,j},\cdots,w_{n,j}]^T$.
Since $\epsilon$ and $\displaystyle \sum_{k=q+1}^p\alpha_kX_k$ are independent, we have for $j=1,\cdots,p-q$ with probability one
as $n\rightarrow \infty$
\begin{eqnarray*}
\frac{1}{n}s_j^2&\rightarrow&
(\sum_{i=1}^n w_{ij}\sum_{k=q+1}^p \frac{x_{ik}{\alpha}_{k}}{\sqrt{n}})^2+\frac{\sigma^2}{n}\\
&\rightarrow&(\sum_{i=1}^n w_{ij}\sum_{k=q+1}^p \frac{x_{ik}{\alpha}_{k}}{\sqrt{n}})^2
\end{eqnarray*}
and 
$${\boldsymbol x}^n:=(\sum_{k=q+1}^p \frac{x_{1,k}\alpha_{k}}{\sqrt{n}},\cdots,\sum_{k=q+1}^p \frac{x_{n,k}\alpha_k}{\sqrt{n}})$$
has a positive constant square norm $||{\boldsymbol x}^\infty||^2$ as $n\rightarrow \infty$ unless $\sum_{k=q+1}^p \alpha_kX_k=0$ with probability one,
which contradicts our assumption.
Since ${\boldsymbol x}^n$ is not orthogonal to the space $<{\boldsymbol w}_1,\cdots,{\boldsymbol w}_{p-q}>$ and $||{\boldsymbol x}^\infty||^2>0$,
from (\ref{eq789}), 
\begin{equation}\label{eq878}
\frac{1}{n}(S_q-S_p)\rightarrow\lim_{n\rightarrow \infty}\sum_{j=q+1}^p({\boldsymbol w}_j^T{\boldsymbol x}^n)^2>0\ ,
\end{equation}
which implies the theorem when $\pi\subset \pi_*$.
Suppose $\pi\not\subset\pi_*$. In the same way, if we notice that (\ref{eq878}) is true even for $q=|\pi\cap \pi_*|$,
so that 
\begin{equation}\label{eq879}
\lim_{n\rightarrow \infty}\frac{1}{n}\{S(\pi\cap \pi_*)-S(\pi_*)\}>0\ .
\end{equation}
Furthermore, if we replace $\pi_*$ by $\pi\cap \pi_*$, from a similar discussion as in Theorem 1, we have
\begin{equation}\label{eq880}
\lim_{n\rightarrow \infty}\frac{1}{n}\{S(\pi)-S(\pi\cap \pi_*)\}=0\ .
\end{equation}
The statements (\ref{eq879})(\ref{eq880}) imply the theorem.
\begin{flushright}(Q. E. D.)\end{flushright}

\section{Proof of the Hannan-Quinn Proposition}

\begin{prop} If $q>p$, with probability one,
\begin{equation}\label{eq89}
\frac{S_p-S_q}{S_p}\leq (q-p)\log\log n
\end{equation}
\end{prop}
Proof: The notation is similar to Proposition 2, and let $p+1\leq j\leq q$.
For $\displaystyle Z_i:=\frac{\sqrt{n}v_{i,j}\epsilon_i}{\sigma}$ with ${\boldsymbol v}_j=[v_{1,j},\cdots,v_{n,j}]^T$, we have
$\displaystyle\sum_{i=1}^n Z_i=\frac{\sqrt{n}{r_j}}{\sigma}$ with expectation zero and variance $\sigma^2$, and 
$E[\sum_{i=1}^nZ_i]=0$, $E[(\sum_{i=1}^nZ_i)^2]=n$. Since
$Z_i$ is independently identically distributed.
$E[Z_i]=0$, $E[Z_i^2]=1$. From the law of iterated logarithms (Stout 1974), we have
$$\frac{\sum_{i=1}^n Z_i}{\sqrt{n\log\log n}}=\frac{{\sqrt{n}}{\boldsymbol v}_{j}^T{\boldsymbol \epsilon}/\sigma}{\sqrt{n\log\log n}}\leq 1\ ,$$
namely,
$$\frac{r_j}{\sigma}\leq \sqrt{\log\log n}$$
with probability one, which means 
$$\frac{S_p-S_q}{S_p/n}\leq (q-p)\log\log n$$
with probability one.
\begin{flushright}(Q. E. D.)\end{flushright}

\begin{teiri}\rm
For $d_n:=2c\log\log n$ ($c>1$), $L(z^n,\pi)>L(z^n,\pi_*)$ with probability one.
\end{teiri}
Proof: 
From Theorem 2, 
the error for $\pi_*\not\subseteq \pi$ is almost surely zero as long as $\displaystyle \frac{d_n}{n}\rightarrow 0$ ($n\rightarrow \infty$),
so that we only need to consider the case $\pi_*\subset \pi$.
However, 
$d_n=2c\log\log n$ with $c>1$  implies the both sides of 
$$\frac{1}{2}\{k(\pi)-k(\pi_*)\}d_n-\frac{1}{4n}[\{k(\pi)-k(\pi_*)\}d_n]^2
\leq n[1-\exp\{-\frac{k(\pi)-k(\pi_*)}{2n}d_n\}]\leq \frac{1}{2}\{k(\pi)-k(\pi_*)\}d_n$$
(see (\ref{eq29})) are at least $(q-p)\log\log n$ with $p=k(\pi_*)$ and $q=k(\pi)$ for large $n$ (Proposition 2), which implies Theorem 3.
\begin{flushright}(Q. E. D.)\end{flushright}

\section{Conclusion}

We proved that the Hannan-Quinn proposition is true for linear regression as well as for auto regression (Hannan-Quinn, 1979)
 and for classification (Suzuki, 2006):
 the minimum rate of $d_n$ satisfying strong consistency is $(2+\epsilon)\log\log n$ for arbitrary $\epsilon>0$.

The future problems contain finding 
 strong consistency conditions that are good for all the cases including linear regression, auto regression,
 and classification. Making clear why the same 
 $d_n=2\log\log n$ is the crucial rate for those problems would be the first step to solve the problem.

\end{document}